\documentclass[12pt]{amsart}
\usepackage{amssymb,amsmath,url,bigstrut,float}
\usepackage{color}
\usepackage[top=1.5in, bottom=1.5in, left=1in, right=1in]{geometry}
\usepackage{algorithmicx}
\usepackage{algorithm}
\usepackage{algpseudocode}
\algrenewcommand\algorithmicrequire{\textbf{Input:}}

\long\def\forget#1\forgotten{}
\newcommand{\nc}{\newcommand}
\nc{\sliding}{\operatorname{sl}}
\nc{\card}[1]{\left|#1\right|}
\nc{\mFunction}[2]{\textsc{#1}(#2)}
\nc{\pos}[1]{#1^+}
\nc{\wt}{\op{wt}}
\nc{\bbZ}{\mathbb{Z}}
\nc{\my}[1]{\textcolor{red}{\sf #1}}
\nc{\set}[2]{\{\, #1 : #2\,\}}
\nc{\tvee}{\mathbin{\tilde\vee}}
\nc{\twedge}{\mathbin{\tilde\wedge}}
\nc{\inv}{^{-1}}
\nc{\BG}[1]{\mathrm{B}_{#1}}
\nc{\Cite}[1]{\textcolor{red}{\textbf{[#1]}}}
\nc{\op}[1]{\operatorname{#1}}
\nc{\suminf}{\op{suminf}}\nc{\sumsup}{\op{sumsup}}\nc{\proj}{\op{proj}}\nc{\Cent}{\op{C}}\nc{\Div}{\op{Div}}
\nc{\pow}{\op{pow}}\nc{\mult}{\op{mult}}\nc{\cl}{\op{cl}}\nc{\Cinf}{\op{C}^\mathrm{inf}}
\renewcommand{\SS}{\op{SS}}\nc{\LSS}{\op{LSS}}\nc{\RSS}{\op{RSS}}\nc{\SSS}{\op{SSS}}
\nc{\LSSS}{\op{LSSS}}\nc{\RSSS}{\op{RSSS}}
\nc{\modtau}{{\mathrm{mod}\, \tau}}
\nc{\rmleft}{\text{left}}\nc{\rmright}{\text{right}}

\newtheorem{thm}{Theorem}[section]
\nc{\bthm}{\begin{thm}} \nc{\ethm}{\end{thm}}
\newtheorem{prop}[thm]{Proposition}
\nc{\bprp}{\begin{prop}} \nc{\eprp}{\end{prop}}
\newtheorem{fact}[thm]{Fact}
\nc{\bfct}{\begin{fact}} \nc{\efct}{\end{fact}}
\newtheorem{prob}[thm]{Problem}
\nc{\bprb}{\begin{prob}} \nc{\eprb}{\end{prob}}
\newtheorem{lem}[thm]{Lemma}
\nc{\blem}{\begin{lem}} \nc{\elem}{\end{lem}}
\newtheorem{claim}[thm]{Claim}
\nc{\bclm}{\begin{claim}} \nc{\eclm}{\end{claim}}
\newtheorem{cor}[thm]{Corollary}
\nc{\bcor}{\begin{cor}} \nc{\ecor}{\end{cor}}
\newtheorem{conj}[thm]{Conjecture}
\nc{\bcnj}{\begin{conj}} \nc{\ecnj}{\end{conj}}
\theoremstyle{definition}
\newtheorem{defn}[thm]{Definition}
\nc{\bdfn}{\begin{defn}} \nc{\edfn}{\end{defn}}
\theoremstyle{remark}
\newtheorem{rem}[thm]{Remark}
\nc{\brem}{\begin{rem}} \nc{\erem}{\end{rem}}
\newtheorem{cnv}[thm]{Convention}
\nc{\bcnv}{\begin{cnv}} \nc{\ecnv}{\end{cnv}}
\newtheorem{exam}[thm]{Example}
\nc{\bexm}{\begin{exam}} \nc{\eexm}{\end{exam}}
\nc{\bpf}{\begin{proof}} \nc{\epf}{\end{proof}}
\nc{\be}{\begin{enumerate}}
\nc{\ee}{\end{enumerate}}
\nc{\bi}{\begin{itemize}}
\nc{\itm}{\item}
\nc{\ei}{\end{itemize}}

\nc{\myalg}[5]{\begin{algorithm}\caption{#2}\label{#1}\begin{algorithmic}
#5\end{algorithmic}\end{algorithm}}

\nc{\ed}{
\subsection*{Acknowledgments}
We thank Patrick Dehornoy and Jean Michel for providing us information about representation theory of Garside groups
and Garside families and, in particular, for discussions that are summarized in
Section~\ref{sec:garsidefamilies}. We thank David Garber for useful discussions.
The research of the first named author is partially
supported by the Minerva Foundation of Germany.

\end{document}
}
\title[Simultaneous Conjugacy in braid groups]{Complete simultaneous conjugacy invariants in Artin's braid groups}

\author[Kalka]{Arkadius Kalka}
\address[Kalka]{Department of Mathematics, Bar-Ilan University, Israel}
\curraddr{Department of Computer Science, Dortmund University of Applied Sciences and Arts, Germany}
\email{arkadius.kalka@fh-dortmund.de}

\author[Tsaban]{Boaz Tsaban}
\address[Tsaban]{Department of Mathematics, Bar-Ilan University, Israel; and
	Faculty of Mathematics, Weizmann Institute of Science, Israel}
\email{tsaban@math.biu.ac.il}

\author[Vinokur]{Gary Vinokur}
\address[Vinokur]{Department of Mathematics, Bar-Ilan University, Israel}
\email{vinokur777@gmail.com}

\begin{document}

\begin{abstract}
We solve the simultaneous conjugacy problem in Artin's braid groups and, more generally, in
Garside groups, by means of
a complete, effectively computable, finite invariant. 
This invariant generalizes the one-dimensional notion of super summit set to arbitrary dimensions. 
One key ingredient in our solution is
the introduction of a provable high-dimensional version of the Birman--Ko--Lee cycling theorem.
The complexity of this solution is a small degree polynomial in the cardinalities
of our generalized super summit sets
and the input parameters.
Computer experiments suggest that the cardinality of this invariant,
for a list of order $N$ independent elements of Artin's braid group $\BG{N}$,
is generically close to~1.
\end{abstract}



\maketitle

\section{Introduction}

In 1911, Dehn formulated three fundamental algorithmic problems concerning groups:
the Word Problem,
the Conjugacy Problem,
and the Group Isomorphism Problem.
The \emph{Word Problem} is that of deciding whether two words in given symmetric
generators of a group represent the same element or, equivalently, whether a word in these generators represents the identity element.
The \emph{Conjugacy Problem} is that of deciding whether two group elements are conjugate.
The \emph{Conjugacy Search Problem} version of this problem is to find, 
given two conjugate group elements, a witness conjugator. 

Throughout, for group elements $g$ and $x$ in $G$, 
we use the notation $g^x:=x\inv gx$.
The \emph{Simultaneous Conjugacy Problem (SCP)} generalizes
the Conjugacy Problem:
$r$-tuples
$(g_1,\dotsc,g_r)$ and $(h_1,\dotsc,h_r)$ of elements of a group $G$ are \emph{conjugate} if there is
an element $x\in G$ such that
\[
(g_1,\dotsc,g_r)^x:=(g_1^x,\dotsc,g_r^x)=(h_1,\dotsc,h_r).
\]
The ($r$-dimensional) SCP is that of deciding whether two $r$-tuples of group elements are conjugate.
The definition of the \emph{Search SCP} is analogous.

An external motivation for studying the SCP comes from cryptography.
The security of a number of cryptographic protocols
reduces to the Search SCP in Artin's braid groups.
In Section~\ref{sec:red}, we provide such reductions for
prominent examples.
In the case where $G$ is a braid group, there
are by now \emph{polynomial-time}, ad-hoc solutions of
the problems that characterize the security of these protocols~\cite{CJ03, LinCent, AlgSpan}.
These solutions do not address the (formally harder) SCP\@.
Moreover, 
for the lack of a general, polynomial-dimension representation theory for Garside groups, the known solutions 
do not generalize to arbitrary Garside groups.
From a \emph{heuristic} point of view,
practically all braid-based cryptographic protocols proposed thus far, 
including ones hitherto not cryptoanalyzed,
are based on the difficulty of the SCP.
In order to understand the potential of braid groups in
cryptography, we must address the full-fledged SCP.

A number of computational problems in braid groups reduce to the SCP. 
For example, Dehornoy's \emph{Shifted Conjugacy Problem}~\cite{Deh06}
reduces to the SCP via a reduction to the Subgroup Conjugacy Problem for the braid group $\BG{N-1}$ in $\BG{N}$~\cite{KLT09}.
More generally, the Subgroup Conjugacy Problem for $\BG{M}$ in $\BG{N}$ ($M<N$) is reducible to the SCP~\cite{GKLT13}.
In a sequel paper~\cite{KTT}, we show that the \emph{Double Coset Problem}
for parabolic subgroups of braid groups reduces to the SCP\@.
In particular, 
the present paper leads to the first solution of the Double Coset Problem.

Our main result is a deterministic, effective solution to the decision and search version of the SCP, in arbitrary Garside groups.
Earlier, Lee and Lee provided a solution in Artin's braid groups~\cite{LL02},
that extends to Garside groups with weighted presentations.
In contrast to the Lee--Lee solution, 
our solution
provides a finite \emph{invariant} of the conjugacy class of an $r$-tuple.
Experimental results, in braid groups,
show a considerable improvement over the earlier solution.
We conclude this paper with open problems and indications
for additional applications.

\section{Reductions of some computational problems to the SCP}\label{sec:red}

In the original instantiations of the problems below, the group $G$ was Artin's braid group. 
The protocols, problems and reductions in this section
apply in arbitrary finitely generated groups.
We assume, for simplicity, 
that each mentioned group is provided in terms of a generating set of cardinality $r$.

The security of the \emph{Braid Diffie--Hellman} protocol~\cite{KL+00}, 
is based on the difficulty of the following problem.%

\bprb
\label{prb:DHCP}
Let $A$ and $B$ be subgroups of a group $G$ with $[A,B]=1$, 
and let $g\in G$ be given.
Given a pair $(g^a, g^b)$, for $a\in A$ and $b\in B$, find $g^{ab}$.
\eprb

Problem~\ref{prb:DHCP} reduces to the Search SCP\@. Indeed,
let $b_1,\dotsc,b_r$ be generators of the subgroup $B$. Find  an element $\tilde a\in G$ such that
\[
(g,b_1,\dotsc,b_r)^{\tilde a}=(g^a,b_1,\dotsc,b_r).
\]
Then $[\tilde a,B]=1$, and
\[
g^{ab}=(g^a)^b=(g^{\tilde a})^b
=g^{\tilde ab}=g^{b\tilde a}=(g^b)^{\tilde a},
\]
which we can compute, having $g^b$ and $\tilde a$.

The \emph{Double Coset} protocol~\cite{CK+01}
is a generalization of the Braid Diffie--Hellman protocol.
Its security is based on the difficulty of the following problem.

\bprb
\label{prb:DCP}
Let $A_1$, $A_2$, $B_1$ and $B_2$ be subgroups of a group $G$ with $[A_1,B_1]=[A_2,B_2]=1$, and let $g\in G$ be given.
Given a pair $(a_1ga_2, b_1gb_2)$, for $a_i\in A_i$ and $b_i\in B_i$ for $i=1,2$,
find $a_1b_1ga_2b_2$.
\eprb

To see that Problem~\ref{prb:DCP} reduces to the Search SCP,
let $b_{i1},\dotsc,b_{ir}$
be generators of the subgroup $B_i$ for $i=1,2$.
For elements $b\in B_1$ of our choice, since $[a_1,B_1]=1$, we know the element $b^{ga_2}=b^{a_1ga_2}$.
Find an element $\tilde a_2\in G$ such that
\[
(b_{11}^g,\dotsc,b_{1r}^g,b_{21},\dotsc,b_{2r})^{\tilde a_2}=
(b_{11}^{ga_2},\dotsc,b_{1r}^{ga_2},b_{21},\dotsc,b_{2r}).
\]
Then $[\tilde a_2,B_2]=1$, and $[(ga_2)(g\tilde a_2)\inv,B_1]=1$.
Compute $\tilde a_1=(a_1ga_2)(g\tilde a_2)\inv=a_1(ga_2)(g\tilde a_2)\inv$. Then $[\tilde a_1,b_1]=1$, and thus
\[
{\tilde a_1}(b_1gb_2){\tilde a_2} =
b_1{\tilde a_1}g{\tilde a_2}b_2 = b_1a_1ga_2b_2= a_1b_1ga_2b_2,
\]
which we can compute from $b_1gb_2$, ${\tilde a_1}$ and ${\tilde a_2}$.

The security of the \emph{Commutator} protocol~\cite{AAG}
is based on the difficulty of the following problem.

\bprb
\label{prb:CKE}
Let $A$ and $B$ be subgroups of a group $G$.
Given $A^b$ and $B^a$, for $a\in A$ and $b\in B$, find the commutator $[a,b]$.
\eprb

The reducibility of Problem~\ref{prb:CKE} to the SCP remains open.
Generically, the centralizer of subgroups of the braid group
is equal to the center of the entire group, which is known.
Problem~\ref{prb:CKE} reduces to the conjunction of the Search SCP and
computing the centralizer of a finite set of elements:
Let $a_1,\dotsc,a_r$, $b_1,\dotsc,b_r$, $c_1,\dotsc,c_r$ and $d_1,\dotsc,d_r$ be generators of the subgroup $A$, $B$ and their centralizers $\Cent(A)$ and
$\Cent(B)$, respectively.
Find elements $\tilde a,\tilde b\in G$ such that
\[
(a_1,\dotsc,a_r,d_1,\dotsc,d_r)^{\tilde b}=(a_1^b,\dotsc,a_r^b,d_1,\dotsc,d_r).
\]
and
\[
(b_1,\dotsc,b_r,c_1,\dotsc,c_r)^{\tilde a}=(b_1^a,\dotsc,b_r^a,c_1,\dotsc,c_r).
\]
Then $a^{\tilde b}=a^b$. Also,
${\tilde a}a\inv\in \Cent(B)$ and $[{\tilde b},\Cent(B)]=1$. In particular,
$[{\tilde b},{\tilde a}a\inv]=1$. Compute
\[
{\tilde a}\inv{\tilde b}\inv{\tilde a}{\tilde b} =
{\tilde a}\inv{\tilde b}\inv({\tilde a}a\inv) a{\tilde b} =
{\tilde a}\inv({\tilde a}a\inv){\tilde b}\inv a{\tilde b} =
a\inv a^{\tilde b} =
a\inv a^b = [a,b].
\]

Finally, the security of the \emph{Centralizer} protocol~\cite{ShUsh06}
is based on the following problem.

\bprb\label{prb:CentKEP}
Let $g,a_1,b_2\in G$, $C\le \Cent(a_1)$, $D\le \Cent(b_2)$, $a_2\in D$ and $b_1\in C$.
Given $g$, $C$, $D$, $a_1ga_2$ and $b_1gb_2$, compute $a_1b_1ga_2b_2$.
\eprb

Problem~\ref{prb:CentKEP} reduces to the conjunction of the Search SCP and
computing the centralizer of a finite set of elements:
Let $c_1,\dotsc,c_r$ and $e_1,\dotsc,e_r$ be generators of the subgroup
$C$ and $\Cent(D)$, respectively.
For each $c\in C$, we can compute $c^{a_1ga_2}=c^{ga_2}$.
Find ${\tilde a_2}\in G$ such that
\[
(c_1^g,\dotsc,c_r^g,e_1,\dotsc,e_r)^{\tilde a_2}=(c_1^{ga_2},\dotsc,c_r^{ga_2},e_1,\dotsc,e_r).
\]
Let ${\tilde a_1}=(a_1ga_2)(g{\tilde a_2})\inv$. Then
${\tilde a_1}=a_1((ga_2)(g{\tilde a_2})\inv)\in \Cent(C)$, and thus
$[{\tilde a_1},b_1]=1$. Also, $[{\tilde a_2},\Cent(D)]=1$ and $b_2\in\Cent(D)$. Thus, $[{\tilde a_2},b_2]=1$. It follows that
\[
{\tilde a_1}(b_1gb_2){\tilde a_2} =
b_1{\tilde a_1}g{\tilde a_2}b_2 =
b_1a_1ga_2b_2 =
a_1b_1ga_2b_2.
\]

\section{Background on Garside groups}

Garside groups~\cite{DP99, Pi01} form a generalization of braid groups where Garside's solution to the conjugacy problem in braid groups
applies.
Many examples of Garside groups are known~\cite{DDGM}. These include, in
addition to Artin's braid groups, all Artin groups of finite type and torus link groups~\cite{Pi03}.

Let $M$ be a monoid. An element 
$a\in M$ is a \emph{left divisor} of an element $b\in M$ 
($a\preceq b$) if $b\in aM$.
An element $a$ is \emph{right divisor of} an element $b$ 
($b\succeq a$) if $b\in Ma$.
It is a \emph{divisor} of $b$ if it is a left or a right divisor of $b$.
$\Div(b)$ is the set of divisors of $b$.
An element $b\in M$ 
is \emph{balanced} if the left and right divisors of $b$ coincide.

A monoid $M$ is \emph{Noetherian} if for each element $a \in M$ there is a
natural number $n$
such that $a$ cannot be expressed as the product of more than $n$ nonidentity elements.
An element $a\ne 1$ in $M$ is an \emph{atom} if $a=bc$ implies $b=1$ or $c=1$.
A set generates a Noetherian monoid $M$ if and only if it includes all atoms of $M$.

Let $M$ be a finitely generated, cancellative Noetherian monoid.
The relations $\preceq$ and $\succeq$ are partial orders,
and every element of $M$ admits only finitely many left and right divisors~\cite{DP99}.
An element $m \in M$ is a \emph{right lcm} of $a$ and $b$
if $a,b\preceq m$, and whenever $a \preceq c$ and $b \preceq c$, we have $m \preceq c$.
The definition of \emph{left lcm} is symmetric, using the relation $\succeq$.
Right and left lcms of pairs $a$ and $b$ are unique, and are denoted $a \vee b$ and $a \tvee b$, respectively.
If $a\vee b$ exists, then there is a unique element $c$ such that $a\vee b=ac$.
This element $c$ is the \emph {right complement (or residue) of $a$ in $b$},  denoted by $a\backslash b$.
We define the \emph{left complement} symmetrically. In particular, we have $a\vee b = a(a\backslash b) = b(b\backslash a)$, and $a \tvee b= (b/a)a = (a/b)b$.

A monoid $M$ is \emph{Gaussian} if it is Noetherian, cancellative, and every
pair of elements $a, b \in M$ admits a right and a left lcm.
Let $M$ be a Gaussian monoid. For every pair of elements $a$ and $b$,
the set of common left divisors of $a$ and $b$ is finite and admits a right lcm,
which is therefore the greatest common left divisor of $a$ and $b$, denoted by $a \wedge b$.
The definition of \emph{right gcd} $a \twedge b$ is symmetric.
A Gaussian monoid $M$ is a lattice with respect to the relations 
$\preceq$ and $\succeq$.
The groups of right fractions and left fractions of the monoid $M$ coincide,
and form the \emph{group of fractions} of $M$. The monoid $M$ embeds in this group.

\bdfn \label{GarsSys}
A \emph{Garside group} is a group $G$ equipped with
a finite subset $S$ and an element $\Delta$ such that
the monoid $G^+$ generated by $S$ is Gaussian with $G$ is its group of fractions,
$\Delta\in G^+$ is a balanced element, and $S=\Div(\Delta)$.

In this case, we say that $G^+$ is a \emph{Garside monoid},
$\Delta $ the \emph{Garside element},
and the elements of $S$ are the \emph{simple elements} of $G$.
\edfn

Let $G$ be a Garside group. There may be several choices of $S$ and $\Delta$ witnessing that, and
we always assume, tacitly, that $S$ and $\Delta$ are fixed in the background.
The set $S$ is closed under the operators $\backslash$, $/$, $\vee$, $\tvee$, $\wedge$ and $\twedge$.
Indeed, $S$ is the closure of the atoms of $G^+$ under 
the operators $\backslash$ and $\vee$.
The functions $\partial\colon a \mapsto a\backslash\Delta$ and $\tilde{\partial }\colon a \mapsto \Delta /a$ map $G^+$ onto $S$, and the restrictions $\partial |_S$, $\tilde{\partial }|_S$ are bijections of $S$ satisfying
$\tilde{\partial }|_S=(\partial |_S)\inv$. In particular, we have 
$\partial^2(a)=\tau (a)$ and $\tilde{\partial }^2=\tau\inv(a)$ for all $a\in S$,
where $\tau$ denotes the inner automorphism of $G$ defined by $a\mapsto \Delta\inv a\Delta$.

The partial orders $\preceq $ and $\succeq$ on $G^+$ extend naturally to partial orders on the whole Garside group $G$.
Define a partial order on $G$ by $a \le b$ if there are $c, c'\in G^+$ such that $b=cac'$.
For elements $a, b \in G$ and $k\in \bbZ$, we have 
$a\preceq \Delta^k \preceq b$ if and only if $b \succeq \Delta^k \succeq a$ if and only if $a \le \Delta^k \le b$.
For $m,n\in\bbZ$, define the \emph{interval}
\[
[m,n] := \set{a\in G}{\Delta^m \le a\le \Delta^n}.
\]
The \emph{infimum} and the \emph{supremum} of an element $a\in G$,
denoted $\inf a$ and $\sup a$, respectively,
are the maximal $m\in\bbZ$ and the minimal $n\in \bbZ$ such that $a\in [m,n]$.
The \emph{canonical length} of $a$, denoted $\cl(a)$, is the difference $\sup a - \inf a$.

Let $G$ be a Garside group. The (left) \emph{normal form} of an element $a\in G^+$
is a unique decomposition $a=s_1\dotsm s_l$
such that $s_i=\Delta \wedge (s_i\dotsm s_l)\in S$ for all $i$, and $s_l\ne 1$.
The length $l$ of this decomposition equals $\sup a$.
The \emph{(left) normal form} of a general element $a$ is obtained by expressing
$a=\Delta^{\inf a}\pos{a}$ for the unique element $\pos{a}\in G^+$,
and decomposing $\pos{a}$ to its (left) normal form.

\section{The basic solution} \label{basic}

Picantin's solution for the Conjugacy Problem in Garside groups~\cite{Pi01}
extends to the SCP by choosing the appropriate coordinate-wise generalization of  the involved notions.
For a natural number $r$, the standard partial order $\le$ of $\bbZ$
extends to a partial order of $\bbZ^r$ coordinate-wise:
for elements $p=(p_1,\dotsc,p_r)$ and $q=(q_1,\dotsc,q_r)$ in $\bbZ^r$, we define
$p\le q$ if $p_i\le q_i$ for all $i=1,\dotsc,r$.
For an $r$-tuple $a=(a_1,\dotsc,a_r)$ of elements of a Garside group $G$,
let
\[
a^G := \set{a^g}{g\in G} = \set{(a_1^g,\dotsc,a_r^g)}{g\in G},
\]
the (simultaneous) \emph{conjugacy class} of $a$.
Define
\begin{align*}
\inf a & := (\inf a_1,\dotsc,\inf a_r);\\
\sup a & := (\sup a_1,\dotsc,\sup a_r).
\end{align*}
For $p,q\in\bbZ^r$ with $p\le q$,
define the following  \emph{interval}:
\begin{align*}
[p,q] 
& = \set{a\in G^r}{p\le\inf a\mbox{ and } \sup a\le q}\\
& = \set{(a_1,\dotsc,a_r)\in G^r}{p_i\le\inf a_i
	\mbox{ and } \sup a_i\le  q_i\mbox{ for all }i=1,\dotsc,r}.
\end{align*}

\blem \label{lem:CinfFin}
Let $G$ be a Garside group, and $r$ be a natural number.
For all tuples $p,q\in \bbZ^r$, the interval $[p,q]$ is finite.
\elem
\bpf
The one-dimensional case ($r=1$) is due to Picantin~\cite{Pi01}.
In the case $r>1$, the interval
\[
[p,q]=[p_1,q_1]\times\dotsb\times [p_r,q_r]
\]
is finite, as a product of finitely many one-dimensional intervals.
\epf

\blem[{\cite[Lemma 2.3]{Pi01}}] \label{L2.3Pic}
Let $G$ be a Garside group. Then
$\Delta \le \alpha \beta$ implies $\Delta\le\alpha (\beta \wedge \Delta )$,
for all $\alpha, \beta \in G^+$.
\elem

Since the set of simple elements $S$ in a Garside group
is finite and the automorphism $\tau \colon S\to S$ a bijection, there is 
a natural number $k$
such that $\tau^k$ is the identity map.
Then the element $\Delta^k$ is in the center of $G$.
Assume that elements $a,c \in G$ are conjugate by an element $x\in G$.
Then, for a large enough natural number $m$, 
the elements $a$ and $c$
are also conjugate by the element $\pos{x}:=\Delta^{mk}x\in G^+$.
It follows that the same holds for tuples $a,c\in G^r$, for an arbitrary 
dimension $r$.

\bthm[Simultaneous Convexity]
\label{simConvTh}
Let $G$ be a Garside group, $p, q\in \bbZ^r$, and $a, c\in [p, q]$.
Assume that $c=a^x=a^{\tilde{x}\inv}$ for elements $x,\tilde{x}\in G^+$.
Let $x_1:=\Delta\wedge x$ and $\tilde{x}_1:=x \twedge\Delta$, the leftmost and rightmost simple factors of $x$, respectively.
Then $a^{x_1},a^{\tilde{x}_1\inv} \in [p,q]$.
\ethm
\bpf
The case $r=1$ is due to Picantin~\cite[Propositions 3.2]{Pi01}. It follows that, for each
$i=1,\dotsc,r$, we have $a_i^{x_1},a_i^{\tilde{x}_1\inv} \in [p_i,q_i]$,
and thus $a^{x_1},a^{\tilde{x}_1\inv} \in [p,q]$.
\epf

The one-dimensional version of the following corollary is due to
Picantin~\cite[Propositions~3.3]{Pi01}.

\bcor\label{PropConvThCoro}
Let $G$ be a Garside group, and $p, q\in \bbZ^r$.
For all conjugate tuples $a, c\in [p,q]$, there are
a natural number $l$, tuples $v_0, v_1, \dotsc, v_l\in [p, q]$,
and simple elements $s_1,\dotsc, s_l\in G$ such that
$v_0=a$, $v_l=c$, and $v_{i-1}^{s_i}=v_i$ for $i=1,\dotsc,l$; schematically:
\[
a\stackrel{s_1}{\longrightarrow }v_1 \stackrel{s_2}{\longrightarrow }v_2
\stackrel{s_3}{\longrightarrow }
\cdots \stackrel{s_{l-1}}{\longrightarrow }v_{l-1} \stackrel{s_l}{\longrightarrow }c.
\]
\ecor
\bpf
There is an element $x\in G^+$ such that $a^x=c$.
The assertion follows by applying the Simultaneous Convexity Theorem~\ref{simConvTh} $\sup x$ times.
\epf

We obtain an extension of Picantin's result~\cite[Corollary 3.4]{Pi01} to the simultaneous setting.

\bthm\label{SCPSolv+}
The SCP in Garside groups is solvable.
\ethm
\bpf
Given tuples $a, c\in G^r$, fix tuples $p,q\in\bbZ^r$ with
$p\le \inf a,\inf c$ and $\sup a,\sup c\le q$.
Then $a, c\in [p,q]$, and the elements $a$ and $c$ 
are conjugate if and only if $c \in a^G\cap [p,q]$.

By Lemma~\ref{lem:CinfFin}, the set $a^G\cap [p,q]$ is finite.
By Corollary~\ref{PropConvThCoro}, this set can be generated by
starting with $a$ and iteratively conjugating with simple elements,
keeping only the conjugates that remain in $[p,q]$, until we obtain no new
elements of $[p,q]$ (Algorithm~\ref{C[p,q]}).

We solve the Search SCP by keeping track of the conjugating elements during
the computation of the set $a^G\cap [p,q]$.
\epf

\myalg{C[p,q]}
{Compute the set $a^G\cap [p,q]$, for $p,q\in\bbZ^r$ and $a\in [p,q]$.}
{}
{}
{\State $W:=\emptyset$; $V:=\{ a \}$
\Repeat
   \State Take $v\in V$
   \ForAll{$s\in S$}
   \State  $u:=(v_1^s,\dotsc,v_r^s)$
   \If{$u\in [p,q]$ and $u\notin V$}
        \State {$V:=V \cup \{u\}$}
   \EndIf
   \EndFor
   \State $W:=W\cup \{v\}$; $V := V\setminus\{v\}$
\Until{$V=\emptyset$}
\State      \Return{$W$}
}

For braid groups, a variation of the
solution presented here was provided by Lee and Lee~\cite{LL02}.
Their solution uses, instead of intervals $[p,q]$, intervals of the form
\[
[p,\infty]  = \set{a\in G^r}{p\le\inf a}.
\]
While these intervals are finite
for braid groups, and more generally for so-called Garside groups with weighted
presentation, they may potentially be infinite in some Garside groups,
in which case the Lee--Lee
solution to the SCP may not terminate in finite time.

The solution presented in this section is infeasible in practice, for two reasons:
The intervals used are typically 
too large, and each step in the algorithm consists of conjugating
by all simple elements. In the braid group $\BG{N}$, there are exponentially (in $N$) many simple elements.
We address these issues in the coming sections.

\section{Simultaneous Cyclic Sliding}

While cycling only affects the infimum of a braid, cyclic sliding~\cite{GG10} affects infimum and supremum.
We identify and establish a high-dimensional generalization of the latter.
This plays a crucial role in our moving to minimal intervals in the next section.

Let $a\in G^r$. For each index $i=1,\dotsc,r$, 
represent the group element $a_i$ in normal form:
\[
a_i=\Delta^{p_i}P_i=\Delta^{p_i}s_1^{(i)} \dotsm s_{l_i}^{(i)}.
\]
Assume that the interval $[p,q]$ is not minimal with respect to $a^G$, that is, 
there exists an element $b\in a^G$ such that
$[p,q] \supsetneqq [\inf b, \sup b ]$.
Consider \emph{target intervals} $[\tilde{p}, \tilde{q}]$ which are proper subintervals of $[p,q]$ such that
$q_i-\tilde{q}_i \le 1$ and $\tilde{p}_i -p_i\le 1$ for all $i=1,\dotsc,r$.
By definition, there are exactly $2r$ target intervals that are maximal with respect to $\le$, namely those where the
tuples $p$ and $\tilde{p}$---(exclusive) or $q$ and $\tilde{q}$---differ in exactly one coordinate.
According to our assumption, there exists (among these $2r$ maximal target intervals) at least one such that
$[\tilde{p}, \tilde{q}] \cap a^G \ne \emptyset $. 
We define simultaneous cyclic sliding with respect to a target interval  $[\tilde{p}, \tilde{q}]$.

\subsection{Simultaneous cycling and decycling}

Let $b \in [\tilde{p}, \tilde{q}] \cap a^G$. 
Fix an element $X \in G^+$ such that $XaX^{-1}=b$. 
For $i=1,\dotsc, r$, write
\begin{equation} \label{bi}
   b_i=\Delta ^{\tilde{p}_i} \bar{b}_i=X\Delta ^{p_i} P_iX^{-1},
\end{equation}
with $\bar{b}_i\succeq 1$. 
The normal form of $b_i$ need not be
$\Delta^{\tilde{p}_i} \bar{b}_i$;
in general, we have $\inf (b_i) \ge \tilde{p}_i$. 
Multiplying Equation~\eqref{bi} on the left by $\Delta ^{-(\tilde{p}_i-1)}$, we have
\[ 
\Delta \bar{b}_i=\Delta ^{-(\tilde{p}_i-1)} X \Delta ^{p_i} P_i X^{-1}=
\tau ^{\tilde{p}_i-1}(X) \Delta ^{p_i-\tilde{p}_i+1} P_i X^{-1} \succeq \Delta.
\]
Since $X\succeq 1$, we have $\tau ^{\tilde{p}_i-1}(X) \Delta ^{p_i-\tilde{p}_i+1} P_i \succeq \Delta$.
For the index $i$ with $\tilde{p}_i=p_i+1$,
we have $\tau ^{\tilde{p}_i-1}(X) P_i \succeq \Delta$ which by 
Lemma~\ref{L2.3Pic} implies that
\[ 
\tau ^{\tilde{p}_i-1}(X) (P_i \wedge \Delta) = \tau ^{\tilde{p}_i-1}(X) s_1^{(i)} \succeq \Delta.
\]
By the invariance of the relation $\succeq $ under right 
multiplication, we have
$ \tau ^{\tilde{p}_i-1}(X) \succeq \Delta (s_1^{(i)})^{-1}=\partial ^{-1} (s_1^{(i)})$.
By invariance of $\succeq$ under $\tau $-automorphism, we have 
$X \succeq \tau ^{-\tilde{p}_i+1}(\partial ^{-1} (s_1^{(i)}))$.
Taking the left lcm for all $i$ with $\tilde{p}_i=p_i+1$, we obtain
\begin{equation} \label{simCyc}
  X \succeq \bigvee _{i: \, \tilde{p}_i=p_i+1}^{\sim } \tau ^{-\tilde{p}_i+1}(\partial ^{-1} (s_1^{(i)})).
\end{equation}

\bdfn
In the above notation, the \emph{simultaneous cycling} operation is the
left conjugation of the tuple $a \in G^r$ by the element in the right hand side of Equation~\eqref{simCyc}.
\edfn

We define simultaneous decycling analogously:
Recall that $q_i=\sup a_i= - \inf a_i^{-1}$. Let the normal form of the element
$a_i$ be $\Delta ^{-q_i} P'_i$.
Since  $\sup b_i \le \tilde{q}_i \le q_i$ and $\inf b_i^{-1}=-\sup b_i$, we have $\inf b_i^{-1}\ge -\tilde{q}_i$.
Thus, for each index $i=1,\dotsc,r$, we can write
\begin{equation} \label{bi-1}
 b_i^{-1}=\Delta ^{-\tilde{q}_i} \bar{b}'_i=Xa_i^{-1}X^{-1}=X\Delta ^{-q_i} P'_iX^{-1}
\end{equation}
for some element $\bar{b}'_i\succeq 1$. 
In general, we have $\inf (b_i^{-1}) \ge -\tilde{q}_i$ and 
$\Delta ^{-\tilde{q}_i} \bar{b}'_i$ need not be the normal form of $b_i^{-1}$.
Multiplying Equation~\eqref{bi-1} on the left by 
$\Delta ^{\tilde{q}_i+1}$, we have
\[
\Delta \bar{b}_i'=\Delta ^{\tilde{q}_i+1} X \Delta ^{-q_i} P'_i X^{-1}=
\tau ^{-\tilde{q}_i-1}(X) \Delta ^{-q_i+\tilde{q}_i+1} P'_i X^{-1} \succeq \Delta.
\]
Since $X\succeq 1$, we have 
$\tau ^{-\tilde{q}_i-1}(X) \Delta ^{-q_i+\tilde{q}_i+1} P'_i \succeq \Delta$.
For $i$ with $\tilde{q}_i=q_i-1$, we have 
$\tau ^{-\tilde{q}_i-1}(X) P'_i = \tau ^{-q_i} (X) P'_i \succeq \Delta$.
By Lemma~\ref{L2.3Pic}, we have
$\tau ^{-q_i}(X) (P'_i \wedge \Delta)\succeq \Delta$.
The normal form of the element $a_i^{-1}$ is related to that of $a_i$.
In particular, we have
\[ 
P'_i \wedge \Delta =\tau ^{-q_i}(\partial (s_{l_i}^{(i)})). 
\]
Thus,
\[ 
\tau ^{-q_i}(X \partial (s_{l_i}^{(i)})) = \tau ^{-q_i}(X (s_{l_i}^{(i)})^{-1} ) \Delta \succeq \Delta.
\]
By invariance of the relation $\succeq$ under right multiplication,
we have
$ \tau ^{-q_i}(X) \succeq \tau^{-q_i}(s_{l_i}^{(i)})$.
Invariance of the relation $\succeq$ under the automorphism 
$\tau$ implies that
$X \succeq s_{l_i}^{(i)}$ for all $i$ with $\tilde{q}_i=q_i-1$.
Finally, we can take the left lcm and obtain
\begin{equation} \label{simDecyc}
 X \succeq \bigvee _{i: \, \tilde{q}_i=q_i-1}^{\sim }  s_{l_i}^{(i)}.
\end{equation}
\bdfn
In the above notation, the \emph{simultaneous decycling} operation is the 
left conjugation of the tuple $a \in G^r$ by the element on the right hand side of Equation~\eqref{simDecyc}.
\edfn

\subsection{Simultaneous cyclic sliding}

\bdfn
In the above notation,
Equations~\eqref{simCyc} and~\eqref{simDecyc} imply that
\begin{equation} \label{simSlid}
 X \succeq \Bigl(\bigvee _{i: \, \tilde{p}_i=p_i+1}^{\sim } \tau ^{-\tilde{p}_i+1}(\partial ^{-1} (s_1^{(i)})) \Bigr)
               \tilde{\vee } \Bigl(\bigvee _{i: \, \tilde{q}_i=q_i-1}^{\sim }  s_{l_i}^{(i)}\Bigr)=:x(0).
\end{equation}
The \emph{simultaneous cyclic sliding} operation 
(with respect to the target interval $[\tilde{p}, \tilde{q}]$)
is the 
left conjugation of the tuple $a \in G^r$ 
by the element $x(0)$.
Let $\sliding(a):=x(0) a x(0)^{-1}$.
\edfn

It is easy to verify that
\[ 
\tilde{p}_i-1\le \inf (\sliding(a_i)) \le \sup (\sliding(a_i)) \le \tilde{q}_i+1
\]
for all $i=1,\dotsc,r$.
In other words, the infimum and supremum of the element $\sliding(a_i)$ 
are at most in distance one outside the target interval.
Indeed, we have implicitly treated the difficult 
cases (where $\tilde{p}_i=p_i+1$ or $\tilde{q}_i=q_i-1$) 
in the derivation above.
The cases where $\tilde{p}_i=p_i$ or $\tilde{q}_i=q_i$ are clear, since conjugation by any simple element (in particular, by $x(0)$) 
can decrease (respectively, increase)
the infimum (respectively, supremum) by at most 1.

For an element $a\in G^+$ in a Garside group $G$, 
let $\|a\|$ be the maximum number of atoms in an expression
of $a$ as a product of atoms.
The following theorem 
generalizes the Birman--Ko--Lee Cycling Theorem~\cite{BKL01}
to dimension $r>1$. It asserts that if moving to a proper subinterval is
possible, then this can be done in at most $\|\Delta\|-1$ steps.
As usual, for $i\in\{1,\dotsc,r\}$ let $e_i\in\bbZ^r$ be the tuple with all coordinates $0$ but
the $i$-th, which is $1$.

\bthm[Simultaneous Cyclic Sliding] \label{SlidThm}
Let $a\in G^r\cap [p,q]$. For each index $i=1,\dotsc,r$ and each
pair $(\tilde p,\tilde q)\in\{(p+e_i,q),(p,q-e_i)\}$ with $a^G\cap [\tilde p,\tilde q]\neq\emptyset$, we have
\[ 
\sliding^{\|\Delta \|-1} (a):=\underbrace{\sliding(\sliding( \cdots \sliding(a)))}_{\|\Delta \|-1 \,\, \text{times}} \in [\tilde{p}, \tilde{q}]. 
\]
\ethm
\bpf
Let $a(0):=a$. 
For $t=0,1,\dotsc,\|\Delta \|-2$, let $a(t+1)=\sliding(a(t))$.
Let $p_i(t):=\inf (a_i(t))$, $q_i(t):=\sup (a_i(t))$, and $l_i(t):=q_i(t)-p_i(t)$. Express the element $a_i(t)$ in normal form:
\[
a_i(t)= \Delta ^{p_i(t)} s_1^{(i)}(t) \dotsm s_{l_i(t)}^{(i)}(t).
\]
Explicitly, $a(t+1)=x(t) a(t) x(t)^{-1}$, where
\[
x(t) := \Bigl( \bigvee _{i: \, \tilde{p}_i=p_i(t)+1}^{\sim } \tau ^{-\tilde{p}_i+1}(\partial ^{-1} (s_1^{(i)}(t))) \Bigr)
               \tilde{\vee } \Bigl( \bigvee _{i: \, \tilde{q}_i=q_i(t)-1}^{\sim }  s_{l_i(t)}^{(i)}(t) \Bigr). 
\]
Setting $X(t):=x(t)x(t-1) \dotsm x(1)x(0)$, we have $a_i(t)=X(t)a_iX(t)^{-1}$.  

Let $m$ and $X=X(m) \succeq 1$ be minimal (with respect to $\succeq$) such that $XaX^{-1}= a(m+1) \in [\tilde{p}, \tilde{q}]$.
Let $\bar{X}(t):=x(m) x(m-1) \dotsm x(t)$, that is, decompose
$X=X(m)=\bar{X}(t)X(t-1)$ for all $t=1, 2, \dotsc, m$. 
For $t=0$, we obtain $\bar{X}(0)=X(m)=X$.
Define 
\[
H(t):=\Delta \wedge \bar{X}(t)
\]
for $t=0,1, \dotsc, m$.
Then $H(t+1) \preceq H(t)$ for all $t=0,1, \dotsc, m-1$. 

By the minimality of the number $m$, we have $H(m)=x(m) \succ 1$. 
By the minimality of the element $X=X(m)$, we have $\inf X=0$, and 
hence $H(0)=\Delta \wedge \bar{X}(0)=\Delta \wedge X \preceq \Delta $. 
In order to prove that
\begin{equation} \label{H}
   1 \prec H(m) \prec \dotsb \prec H(1) \prec H(0) \prec \Delta
\end{equation}
it suffices to show that $H(t) \ne H(t+1)$ for all $t$.
Then Equation~\eqref{H} implies that
\[ 
0 < \|H(m)\| < \dotsb < \|H(1)\| < \|H(0)\| < \|\Delta \|,
\]
and thus $m+1$, is bounded below $\|\Delta \|-1$, which completes the proof.  

We prove the inequality $H(1) \ne H(0)$;
the proof of the inequality
$H(t+1) \ne H(t)$ for $t=1,\dotsc, m-1$ is similar.
Let $L:=\cl(X)=\sup(X) \le m+1$.
Express in normal form $X=B_L \dotsm B_2 B_1$. 
Then $H(0)=\bar{X}(0) \wedge \Delta =X \wedge \Delta =B_L$. 
Assume, towards a contradiction, that $H(1)=H(0)=B_L$, that is, $\bar{X}(1) \wedge \Delta =B_L$.
Write $\bar{X}(1)=x(m) \dotsm x(1)=B_L R_1$ for some $R_1 \in G^+$. 
Then $X=B_L R_1 x(0)$ and $B_{L-1} \dotsm B_1=R_1 x(0)$. 
We prove that
\begin{align}
 \label{Ba} (B_{L-1} \dotsm B_1) a_i &\succeq \Delta ^{\tilde{p}_i}, \text{and} \\
 \label{Ba-1} (B_{L-1} \dotsm B_1) a_i^{-1} &\succeq \Delta ^{-\tilde{q}_i}
\end{align}
For all $i=1,\dotsc,r$.

Proof of Equation~\eqref{Ba}. First, let $i$ be an index with $\tilde{p}_i=p_i+1$.
According to Equation~\eqref{simCyc}, we have $x(0) \succeq \tau ^{-p_i}(\partial ^{-1} (s_1^{(i)}))$.
Write $x(0)=r_i \tau ^{-p_i}(\partial ^{-1} (s_1^{(i)}))$. Then
\begin{align*}
 (B_{L-1} \dotsm B_1) a_i &= R_1 x(0) \Delta ^{p_i} P_i= R_1 r_i \tau ^{-p_i}(\partial ^{-1} (s_1^{(i)}))\Delta ^{p_i}P_i \\
&= R_1 r_i \Delta ^{p_i}  \underbrace{\partial ^{-1} (s_1^{(i)}) s_1^{(i)}}_{\Delta } s_2^{(i)} \dotsm s_{l_i}^{(i)} \succeq \Delta ^{p_i+1}
=\Delta ^{\tilde{p}_i}.
\end{align*}
The case where $\tilde{p}_i=p_i$ is simpler:
\[ 
(B_{L-1} \dotsm B_1) a_i =R_1 x(0) \Delta ^{p_i} P_i \succeq \Delta ^{p_i}=\Delta ^{\tilde{p}_i}. 
\]

Proof of Equation~\eqref{Ba-1}: Let $i$ be an index such that $\tilde{q}_i=q_i-1$.
According to Equation~\eqref{simDecyc}, we have $x(0) \succeq s_{l_i}^{(i)}$.
Write $x(0)=r'_i s_{l_i}^{(i)}$. Since $a_i^{-1}=\Delta ^{-q_i} \bar{P}_i$ 
and $\bar{P}_i \wedge \Delta =\tau ^{-q_i}(\partial (s_{l_i}^{(i)}))$), we can write
$\bar{P}_i =\tau ^{-q_i}(\partial (s_{l_i}^{(i)})) \bar{P}'_i$.
Thus,
\begin{align*}
 (B_{L-1} \dotsm B_1) a_i^{-1} &= R_1 x(0) \Delta ^{-q_i} \bar{P}_i= R_1 r'_i s_{l_i}^{(i)} \Delta ^{-q_i} \bar{P}_i \\
&= R_1 r'_i \Delta ^{-q_i}  \underbrace{\tau ^{-q_i} (s_{l_i}^{(i)}) \tau ^{-q_i}(\partial (s_{l_i}^{(i)}))}_{\Delta } \bar{P}'_i \succeq \Delta ^{-q_i+1}
=\Delta ^{-\tilde{q}_i}.
\end{align*}
The case where $\tilde{q}_i=q_i$ is simpler:
\[ 
(B_{L-1} \dotsm B_1) a_i^{-1} =R_1 x(0) \Delta ^{-q_i} \bar{P}_i \succeq \Delta ^{-q_i}=\Delta ^{-\tilde{q}_i}.
\]
Next, we prove that:
\begin{align}
 \label{infBa}  \inf ((B_{L-1} \dotsm B_1) a_i (B_{L-1}\dotsm B_1)^{-1}) &\ge  \tilde{p}_i \text{ and } \\
 \label{supBa} \sup ((B_{L-1} \dotsm B_1) a_i (B_{L-1}\dotsm B_1)^{-1}) &\le  \tilde{q}_i 
\end{align}
For all $i=1,\dotsc,r$.

Proof of Equation~\eqref{infBa}: Define, for $k=0,1, \dotsc, L$ and $i=1,\dotsc,r$,
\[ \alpha _i(k):=\inf [(B_{L-1} \dotsm B_1) a_i \partial (B_1) \tau (\partial (B_2)) \dotsm \tau ^{k-1}(\partial (B_k)) ]. \]
Note that $\partial (B_1) \tau (\partial (B_2)) \dotsm \tau ^{k-1}(\partial (B_k))$ is the normal form of
$(B_k \dotsm B_1)^{-1} \Delta ^k$. 
By definition, for $i=1,\dotsc,r$ we have 
\begin{align*}
\alpha _i(0) &= \inf [(B_{L-1} \dotsm B_1) a_i] \ge \tilde{p}_i, \\
\alpha _i(L) &= \inf [\underbrace{(B_{L-1} \dotsm B_1)}_{B_L^{-1}X} a_i
                \underbrace{\partial (B_1)\dotsm \tau ^{L-1}(\partial (B_L))}_{X^{-1}\Delta ^L}] \\
						 &= \inf [B_L^{-1}(Xa_iX^{-1})\Delta ^L] \ge -1 + \tilde{p}_i +L, \quad \text{and} \\
\alpha _i(k) &\le \alpha _i(k+1)						
\end{align*}
for all $k=0,1,\dotsc L-1$.
By the lemma, we have
$\alpha _i(k+1) \le \alpha _i(k)+1$ for all $k$ and $i$.
Since $\partial (B_1) \tau (\partial (B_2)) \dotsm \tau ^{k-1}(\partial (B_k))$ is in normal form, we have
\[ 
\alpha _i(k)=\alpha _i(k+1) \quad \Rightarrow \quad \alpha _i(k)= \dotsb =\alpha _i(L-1)=\alpha _i(L). 
\]
Assume that there is an index $i$ such that $\alpha _i(L-1) \le \tilde{p}_i+L-2$. 
Since $\alpha _i(0) \ge \tilde{p}_i$,
there is a natural number $k$ such that $\alpha _i(k)=\alpha _i(k+1)$. Hence
$\alpha _i(k)= \dotsb =\alpha _i(L-1)=\alpha _i(L) \le \tilde{p}_i+L-2$, in contradiction to the inequality $\alpha _i(L)\ge \tilde{p}_i +L-1$.
Thus, $\alpha _i(L-1) \ge \tilde{p}_i +L-1$ for all $i=1,\dotsc, r$, that is,
\[ 
\inf[ (B_{L-1} \dotsm B_1) a_i
\underbrace{\partial (B_1) \dotsm \tau ^{L-2}(\partial (B_{L-1}))}_{(B_{L-1}\dotsm B_1)^{-1}\Delta ^{L-1}} ] \ge \tilde{p}_i +L-1,
\]
or, equivalently,
$\inf[(B_{L-1} \dotsm B_1) a_i (B_{L-1}\dotsm B_1)^{-1} ] \ge \tilde{p}_i$ 
for all $i$.

Proof of Equation~\eqref{supBa}:
Analogously, define
\[ 
\alpha _i^-(k):=\inf [(B_{L-1} \dotsm B_1) a_i^{-1} \partial (B_1) \tau (\partial (B_2)) \dotsm \tau ^{k-1}(\partial (B_k)) ]. 
\]
The element $\alpha _i^-(k)$ has the same properties as 
$\alpha _i(k)$, expect that we have to replace $\tilde{p}_i$ 
by $-\tilde{q}_i$.
In particular, we have $\alpha _i^-(0) \ge -\tilde{q}_i$ and 
$\alpha _i^-(L) \ge -1-\tilde{q}_i+L$. 
The proof proceeds as in the previous case, and we obtain
\[ 
\inf [(B_{L-1} \dotsm B_1) a_i^{-1} (B_{L-1}\dotsm B_1)^{-1}]=-\sup [(B_{L-1} \dotsm B_1) a_i (B_{L-1}\dotsm B_1)^{-1}] \ge -\tilde{q}_i.
\]

Equations~\eqref{infBa} and~\eqref{supBa} assert that the tuple 
$(B_{L-1} \dotsm B_1)a(B_{L-1} \dotsm B_1)^{-1}$ 
lies in the target interval
$[\tilde{p}, \tilde{q}]$,
in contradiction to the minimality of the element
$X=B_L B_{L-1} \dotsm B_1$. 
We conclude that the assumption $H(0)=H(1)$ is false.
\epf

\section{Moving to minimal intervals}

Let $G$ be a Garside group and $a,c\in G^r$.
The solution to the SCP  for $a$ and $c$,
described in section~\ref{basic}, is by choosing some interval $[p,q]$
containing $a$ and $c$ and computing the set
$a^G\cap [p,q]$. Increasing $p$ or decreasing $q$---lexicographically---may
reduce the cardinality of this set considerably.

Algorithm~\ref{LLAlgoInt} conjugates an $r$-tuple into a prescribed interval,
assuming that this is possible.
For braid groups, with $[\inf c,\infty]$ (where $c$ is
conjugate to $a$) instead of $[p,q]$ (and $M=\infty$),
this algorithm is similar to that of 
Lee and Lee~\cite{LL02}. 
Here, for example, we can take the smaller interval $[\inf c,\sup c]$.
We note that a simultaneous cycling theorem was not established 
for the operation used by Lee and Lee; they did not conjugate by the 
$\tilde{\vee }$-join of all cyclings for components with infimum outside the target interval.

\myalg{LLAlgoInt}
{Given tuples $a \in G^r$ and $p,q\in\bbZ^r$, find an element $y\in G^+$ such that
$a^{y\inv}\in [p,q]$. Uses input parameter $M\in\mathbb{N}\cup\{\infty\}$.}
{}
{}
{\Function{ConjugateToInterval}{$a,p,q,M$}
\State $y:=1$; $c:=a$; $i:=0$
\While {$c\notin [p,q]$ and $i<M$}
   \State $h:=1$
   \For{$k:=1$ to $r$}
     \If{$\inf c_k < p_k$}
     \State $h:=h \tvee \Delta \tau^{-\inf c_k}(\Delta\wedge(\Delta^{-\inf c_k}c_k))\inv$
     \EndIf
     \If{$q_k < \sup c_k$}
     \State Bring $c_k$ in normal form $\Delta ^{\inf c_k} s_1 \dotsm s_{\cl(c_k)}$
     \State $h:=h \tvee s_{\cl(c_k)}$
     \EndIf
   \EndFor
   \State $y:=hy$, $c:=c^{h\inv}$
   \State $i:=i+1$
\EndWhile
\State \Return{$y$}\Comment{Successful if and only if $a^{y\inv}\in [p,q]$.}
\EndFunction
}

By the Simultaneous Cyclic Sliding Theorem (Theorem~\ref{SlidThm}),
we have the following performance guarantee.

\bcor
Assume that $a\in G^r\cap [p,q]$, $i\in\{1,\dotsc,r\}$, and 
$(\tilde{p},\tilde{q})\in\{(p+e_i,q),(p,q-e_i)\}$ is a pair 
with $a^G\cap [\tilde {p},\tilde{q}]\neq\emptyset$. Let
$y:= {}$\mFunction{ConjugateToInterval}{$a,p',q',\|\Delta\| -1$}.
Then $a^{y\inv}\in [\tilde{p},\tilde{q}]$.
\ecor

\bdfn  \label{MinInt}
Let $G$ be a Garside group, $r$ a natural number, $a\in G^r$,
and $p,q\in \bbZ^r$ be tuples with $p\le q$.
The interval $[p,q]$
is \emph{minimal} for the conjugacy class $a^G$ if $[p,q]$ intersects $a^G$,
but no proper subinterval of $[p,q]$ intersects $a^G$.
\edfn

Consider the one-dimensional case. For a Garside group
$G$ and an element $a\in G$, the \emph{summit infimum} and 
\emph{summit supremum} of $a$
are the maximal infimum and minimal supremum, respectively, of an element
of $a^G$. In this one-dimensional case,
the interval $[\suminf(a), \sumsup(a)]$ is the only minimal interval with respect to $a$.
Garside's \emph{Summit Set}~\cite{Ga69} and Elrifai--Morton's 
\emph{Super Summit Set}~\cite{EM94} of $a$ are the sets
\begin{align*}
\SS(a) & = a^G\cap [\suminf a,\infty]; \\
\SSS(a) & = a^G\cap [\suminf a,\sumsup a],
\end{align*}
respectively.
These sets  are complete conjugacy invariants.
In higher dimensions, there are in general more than one minimal interval for a
conjugacy class. Any canonical choice among them would provide a complete
conjugacy invariant. 
We provide two variations of a complete invariant for simultaneous
conjugacy classes. Since we use minimal intervals, these invariants generalize
the classic Super Summit Sets to general dimension.

\bdfn\label{LSSS}
Let $G$ be a Garside group and  $a \in G^r$.
The \emph{lexicographically minimal interval} for the
conjugacy class $a^G$ is the unique interval $[p,q]$ with the following properties:
\be
\item $p_1$ and $q_1$ are the summit infimum and summit supremum of $a_1$,
respectively.
\item For $i=2,3,\dotsc,r$, \emph{in this order}: $p_i$ is maximal and $q_i$ is minimal (given $p_i$) with
$(a_1,\dotsc,a_i)^G\cap [(p_1,\dotsc,p_i),(q_1,\dotsc,q_i)]\ne\emptyset$.
\ee
The \emph{Lexicographic Super Summit Set} of $a$, $\LSSS(a)$, is the intersection
of $a^G$ with its lexicographically minimal interval.

The \emph{lexicographically$'$ minimal interval} for the
conjugacy class $a^G$ is the unique interval $[p,q]$ with the following properties:
\be
\item For $i=1,\dotsc,r$: $p_i$ is maximal with
$(a_1,\dotsc,a_i)^G\cap [(p_1,\dotsc,p_i),\infty]\ne\emptyset$.
\item For $i=1,\dotsc,r$: $q_i$ is minimal with
$(a_1,\dotsc,a_i)^G\cap [(p_1,\dotsc,p_i),(q_1,\dotsc,q_i)]\ne\emptyset$.
\ee
The \emph{Lexicographic$'$ Super Summit Set} of $a$, $\LSSS'(a)$, is the intersection
of $a^G$ with its lexicographically$'$ minimal interval.
\edfn
If $[p,q]$ is the lexicographically$'$ minimal interval for a conjugacy
class $a^G$, then the set $\LSS(a):=a^G\cap [p,\infty]$ is also a complete invariant for $a^G$,
but it is, in general, larger than $\LSSS'(a)$.

Given a lexicographically minimal interval $[p,q]$ for an element $a\in G^r$,
we can conjugate the element $a$ into $\LSSS(a)$ using
\mFunction{ConjugateToInterval}{$a,p,q,\|\Delta\| -1$}, and then compute the entire
$\LSSS(a)$ using Algorithm~\ref{C[p,q]}. 
Algorithm~\ref{LMI} computes the lexicographically minimal interval
for an element $a$.
Analogous assertions hold for $\LSSS'(a)$.

\myalg{LMI}
{Compute the lexicographically minimal interval $[p,q]$ for an element $a\in G^r$}
{}
{}
{\For{$i=1,\dotsc,r$}
        \State $p_i:= \inf a_i$.
        \Repeat
        \State $y := {}$\mFunction{ConjugateToInterval}
        {$(a_1,\dotsc,a_i),(p_1,\dotsc,p_i+1),(q_1,\dotsc,q_{i-1},\infty),\|\Delta\| -1$}
        \If{$(a_1,\dotsc,a_i)^{y\inv}\in [(p_1,\dotsc,p_i+1),(q_1,\dotsc,q_{i-1},\infty)]$}
                \State $p_i := p_i+1$
                \State $a:=a^{y\inv}$
                \State \textsf{flag${} := {}$true}
        \Else
                \State \textsf{flag${} := {}$false}
        \EndIf
        \Until \textsf{flag${} = {}$false}
        \State $q_i:= \sup a_i$.
        \Repeat
        \State $y := {}$\mFunction{ConjugateToInterval}
        {$(a_1,\dotsc,a_i),(p_1,\dotsc,p_i),(q_1,\dotsc,q_i-1),\|\Delta\| -1$}
        \If{$(a_1,\dotsc,a_i)^{y\inv}\in [(p_1,\dotsc,p_i+1),(q_1,\dotsc,q_i-1)]$}
                \State $q_i := q_i-1$
                \State $a:=a^{y\inv}$
                \State \textsf{flag${} := {}$true}
        \Else
                \State \textsf{flag${} := {}$false}
        \EndIf
        \Until \textsf{flag${} = {}$false}
\EndFor
\State \Return{$p,q$}
}

The following proposition summarizes the relations among the introduced
invariants.

\bprp\label{InvSubsets}
Let $G$ be a Garside group and $i\in \{1,\dotsc,r\}$.
Let $\proj_i$ denote the projection on the first $i$ coordinates.
For each tuple $a \in G^r$, the following relations hold:
\begin{align*}
\LSS((a_1, \dotsc, a_i)) &= \proj_i(\LSS(a)), \\
\LSSS((a_1, \dotsc, a_i)) &= \proj_i(\LSSS(a)), \\
\LSSS'((a_1, \dotsc,a_i)) &\subseteq \proj_i(\LSSS'(a)), \\
\LSSS'(a) &\subseteq \LSS(a).
\end{align*}
In particular, we have $\SS(a_1)=\proj_1(\LSS(a))$ and
\[ 
\SSS(a_1)=\proj_1(\LSSS(a)) \subseteq \proj_1(\LSSS'(a)).\qed  
\]
\eprp

Our invariants are computable in finite time:
detecting the lexicographically minimal interval $[p,q]$ of an element of this set, and then computing $a^G\cap
[p,q]$.
Also, note that finite invariants of conjugacy
classes imply canonical representatives: The lexicographically minimal
element of the invariant. However, the computational complexity of computing
such a canonical representative remains proportional to the cardinality of the initial invariant.

This completes our treatment of interval minimization. We next
address the second and last problem: Removing the need to conjugate by all
simple elements in each step of our algorithms. This will be done by extending
the method of minimal simple elements to our situation.

\section{Minimal simple elements}

We apply the technique of minimal simple elements, introduced by Gonz\'alez--Meneses and Franco~\cite{FG03}, 
in order to make the computation of the sets $a^G\cap [p,q]$ more efficient.
The propositions and algorithms in this section are natural
generalizations of earlier algorithms~\cite{FG03, Go05, KLT10}.

\bprp \label{charPv}
Let $G$ be a Garside group.
Let $v\in G^r$ be a tuple with $v\in [p,q]$, for $p,q\in \bbZ^r$.
For $i=1,\dotsc,r$,
express $v_i=\Delta^{p_i}w_i=z_i\Delta^{p_i}$ with $w_i, z_i\in G^+$ and
$(w_i\inv\Delta^{q_i-p_i})$, $(\Delta^{q_i-p_i}z_i\inv) \in G^+$. Let $s\in S$.
Then:
\be
\item
$v^s\in [p,q]$ if and only if
\[
\tau^{p_i}(s) \preceq w_is \mbox{ and }\tau^{-q_i}(s)\preceq (\Delta^{q_i-p_i}z_i\inv)s
\]
for all $i=1,\dotsc,r$.
\item
$v^{s\inv}\in [p,q]$ if and only if
\[
sz_i \succeq \tau^{-p_i}(s) \mbox{ and } s(w_i\inv\Delta^{q_i-p_i}) \succeq \tau^{q_i}(s)
\]
for all $i=1,\dotsc,r$.
\ee
\eprp
\bpf
(1) For $i=1,\dotsc,r$, we have 
$\Delta^{p_i} \le v_i^s \le \Delta^{q_i}$ if and only if
$\Delta^{p_i} \preceq v_i^s \preceq \Delta^{q_i}$.
Since $\preceq $ is invariant under left multiplication, we have
$\Delta^{p_i} \preceq s\inv\Delta^{p_i} w_is$ if and only if
$\tau^{p_i}(s) \preceq w_is$,
and $ s\inv z_i\Delta^{p_i} s \preceq \Delta^{q_i}$ is equivalent to
$s\preceq \Delta^{-p_i}z_i\inv s\Delta^{q_i}=\tau^{q_i}(\Delta^{q_i-p_i}z_i\inv s)$.
By invariance of $\preceq $ under the automorphism $\tau$, we have 
$\tau^{-q_i}(s)\preceq (\Delta^{q_i-p_i}z_i\inv)s$.

(2) Here, for $i=1,\dotsc,r$, we use that $\Delta^{p_i} \le v_i^{s\inv} \le \Delta^{q_i}$
if and only if $\Delta^{q_i} \succeq v_i^{s\inv} \succeq \Delta^{p_i}$.
Since $\succeq $ is invariant under right multiplication, we conclude that
$sz_i\Delta^{p_i}s\inv \succeq \Delta^{p_i}$ if and only if
$sz_i \succeq \tau^{-p_i}(s)$,
and $\Delta^{q_i} \succeq s\Delta^{p_i}w_is\inv$ is equivalent to
$\Delta^{q_i}sw_i\inv\Delta^{-p_i}=\tau^{-q_i}(sw_i\inv\Delta^{q_i-p_i}) \succeq s$.
By invariance of $\succeq $ under $\tau$, we have 
$s(w_i\inv\Delta^{q_i-p_i}) \succeq \tau^{q_i}(s)$.
\epf


\blem \label{Lem4.2}
Assume that  a set $A\subseteq S$ is closed under the operation
$\wedge$ (respectively, $\twedge$).
Let $x \in A$.
If the set $\set{ s\in A }{ x \preceq s }$ (respectively, $ \set{ s\in A }{ s \succeq x }$)
is nonempty, then it has a unique minimal element with respect to the relation $\preceq $
(respectively, $\succeq $).
\elem
\bpf
Every interval $\set{ x\in G^+ }{ a \preceq x\preceq b }$
in the poset $(G^+,\preceq )$ is closed under the operations $\wedge $ and $\vee $.
The intersection of sets closed under $\wedge$ and $\vee$ is also closed
under these operations.
Uniqueness follows.
\epf

\bdfn \label{Smin}
Let $G$ be a Garside group, $p,q\in\bbZ^r$, and $v\in [p,q]$.
The set $S_\rmright^{[p,q]}(v)$ consists of all $\preceq$-minimal 
elements $s\in S$
such that $v^s\in [p,q]$. Similarly, the set $S_\rmleft^{[p,q]}(v)$
consists of all $\succeq$-minimal elements $s\in S$
such that $v^{s\inv}\in [p,q]$.
\edfn

Analogously to the proof of Gonz\'ales--Meneses~\cite[Proposition 2.2]{Go05}, 
we prove the following result.

\bthm \label{Theorem}
Let $G$ be a Garside group, and $a, c\in [p,q]\subseteq G^r$.
The following assertions are equivalent:
\be
\item The tuples $a$ and $c$ are conjugate.
\item There exist a natural number $l$, elements
$\tilde{w}_1,\dotsc, \tilde{w}_{l-1}\in a^G\cap[p,q]$, and elements
$\tilde{s}_i \in S_\rmright^{[p,q]}(\tilde{w}_{i})$, for $i=1,\dotsc,l$,
such that
\[ a\stackrel{\tilde{s}_1}{\longleftarrow }\tilde{w}_1
\stackrel{\tilde{s}_2}{\longleftarrow }
\cdots
\stackrel{\tilde{s}_{l-1}}{\longleftarrow }
\tilde{w}_{l-1}
\stackrel{\tilde{s}_l}{\longleftarrow }c. \]
\item There exist  a natural number $l$, 
elements $w_1,\dotsc, w_{l-1}\in a^G\cap [p,q]$,
and elements $s_i \in S_\rmleft^{[p,q]}(w_i)$, for $i=1,\dotsc,l$, such that
\[
a
\stackrel{s_1}{\longrightarrow }
w_1
\stackrel{s_2}{\longrightarrow }
\cdots
\stackrel{s_{l-1}}{\longrightarrow }
w_{l-1}
\stackrel{s_l}{\longrightarrow }
c.\qed
\]
\ee
\ethm

The set $\set{s\in S}{v^s\in [p,q]}$ is closed under $\wedge$.
It follows~\cite[Corollary 4.3]{FG03} that
for each $s\in S_\rmright^{[p,q]}(v)$ there is an atom $x\in S$
such that $s$ is the unique $\preceq$-minimal element of the set
$\set{a\in S}{x\preceq a\mbox{ and }v^a\in [p,q]}$.
Similarly, each element of
$S_\rmleft^{[p,q]}(v)$ is the unique $\succeq$-minimal element of the set
$\set{a\in S}{a\succeq x\mbox{ and }v^a\in [p,q]}$
for some atom $x\in S$.
It follows that, in Algorithm~\ref{C[p,q]}, the computation of $a^G\cap [p,q]$
can be done with $S$ replaced by $S_\rmright^{[p,q]}(v)$, a set not larger
than the number of \emph{atoms} in $S$.
For example, the Artin groups of type $A_n$, $B_n$, and $D_n$, respectively,
have $(n+1)!$, $n!2^n$, and $n!2^{n-1}$ simple elements, but only $n$ atoms.


Algorithms~\ref{rxright},~\ref{rxleft} and~\ref{Alg3},
build on earlier algorithms~\cite{FG03, Go05,KLT10}.

\myalg{rxright}
{Compute the minimal element in the set 
$\set{a\in S}{x\preceq a\mbox{ and }v^a\in [p,q]}$,
for an atom $x$.}
{Atom $x$ in $S$, and an element $v\in [p,q]\subseteq G^r$.}
{The unique minimal element $r_x$ in $S_{\mathcal{P}_\rmright^{[p,q]}(v)} \cap \mult_{\preceq}(x)$.}
{\State Express $v_i=\Delta^{p_i} w_i=z_i\Delta^{p_i}$
\State $s:=x$
\While{There is an index $k\in \{1,\dotsc, r\}$ with $\tau^{p_k}(s) \npreceq w_ks$
or
$\tau^{-q_k}(s)\npreceq (\Delta^{q_k-p_k}z_k\inv)s$}
\State   Choose such an index $k$
\State   Compute $s_1\in S$ such that $\tau^{p_k}(s)\vee w_ks=w_kss_1$
\State   $w'_k:=\Delta^{q_k-p_k}z_k\inv$
\State   Compute $s_2\in S$ such that $\tau^{-q_k}(s)\vee w'_ks=w'_kss_2$
\State   $s':=s_1\vee s_2$
\State   $s:=ss'$
\State $r_x:=s$
\EndWhile
\State \Return $r_x$
}

\bprp \label{rxrAlgworks}
Let $G$ be a Garside group, $p,q\in \bbZ^r$, $v\in [p,q] \subseteq G^r$, and $x\in S$.
Algorithm~\ref{rxright} terminates and provides the correct output.
\eprp
\bpf
By Proposition~\ref{charPv} (1), we need to find
the smallest element $r_x$ such that for all $i=1,\dotsc,r$, $x\preceq r_x$, $\tau^{p_i}(s) \preceq w_is$ and
$\tau^{-q_i}(s)\preceq w'_is$ with $w'_i=\Delta^{q_i-p_i}z_i\inv$.
We take a simple element $s$ such that $x\preceq s \preceq r_x$, initializing with $s:=x$.
Then, for $k$ such that $\tau^{p_k}(s) \npreceq w_ks$ or
$\tau^{-q_k}(s)\npreceq (\Delta^{q_k-p_k}z_k\inv)s$, we compute $s_1, s_2\in G$ such that
$\tau^{p_k}(s)\vee w_ks=w_kss_1$ and $\tau^{-q_k}(s)\vee w'_ks=w'_kss_2$.
If $s':=s_1\vee s_2=1$, then $s_1=s_2=1$, and we have $\tau^p(s)\preceq ws$ and
$\tau^{-q}(s)\preceq w's$, in contradiction to the choice of $k$.
Thus, $s'\ne 1$. $s\preceq r_x$ implies that
$\tau^{p_k}(s)\preceq \tau^{p_k}(r_x)\preceq w_kr_x$, and by left-invariance of $\preceq$,
it also implies that  $w_ks\preceq w_kr_x$, that is,
$\tau^{p_k}(s),w_ks\preceq w_kr_x$. By the definition of right lcm, we have
that
$\tau^{p_k}(s)\vee w_ks=w_kss_1\preceq w_kr_x$, and therefore $ss_1\preceq r_x$.

Furthermore,  $s\preceq r_x$ implies that  $\tau^{-q_k}(s)\preceq \tau^{-q_k}(r_x)\preceq w'_kr_x$, and by
left-invariance of the relation
$\preceq $ it also implies $w'_ks\preceq w'_kr_x$, that is,
$\tau^{-q_k}(s),w'_ks\preceq w'_kr_x$.
By the definition of right lcm, 
$\tau^{-q_k}(s)\vee w'_ks=w'_kss_2\preceq w'_kr_x$, and therefore
$ss_2\preceq r_x$.

From $ss_1\preceq r_x$ and $ss_2\preceq r_x$ we conclude that  $ss_1\vee ss_2=s(s_1\vee s_2)=ss'\preceq r_x$.

So, if $s$ is not equal to $r_x$, then Algorithm~\ref{rxright} gives an element $s'\ne 1$ such that
$s\preceq ss' \preceq r_x$,
and it starts again checking whether $ss'= r_x$. Since the number of left divisors
of $r_x$ is finite, this process must stop. Therefore, Algorithm~\ref{rxright} finds
the requested minimal element in finite
time.
\epf

\myalg{rxleft}
{Compute the minimal element of the set $\set{a\in S}{a\succeq x\mbox{ and }v^a\in [p,q]}$,
for an atom $x\in S$}
{Atom $x$ in $S$, $p,q\in \bbZ^r$, and an element $v\in [p,q]\subseteq G^r$.}
{The unique minimal element $r_x$ in $S_{\mathcal{P}_\rmleft^{[p,q]}(v)} \cap \mult_{\succeq}(x)$.}
{\State Express $v_i=\Delta^{p_i} w_i=z_i\Delta^{p_i}$
\State  $s:=x$
\While{There is $k$ with $sz_k \nsucceq \tau^{-p_k}(s)$
and
$s(w_k\inv\Delta^{q_k-p_k}) \succeq \tau^{q_k}(s)$}
\State Compute $s_1\in S$ such that $\tau^{-p_k}(s)\tilde{\vee } sz_k=s_1sz_k$
\State $w'_k:=w_k\inv\Delta^{q_k-p_k}$
\State Compute $s_2\in S$ such that $\tau^{q_k}(s)\tilde{\vee } sw'_k=s_2sw'_k$
\State $s':=s_1 \tilde{\vee } s_2$
\State $s:=s's$
\EndWhile
\State $r_x:=s$
\State \Return $r_x$
}

The sets of minimal simple elements $S_\rmright^{[p,q]}(v)$ and $S_\rmleft^{[p,q]}(v)$
can be computed by comparing the elements $r_x$, for all atoms $x$
of $G^+$, and keeping the minimal ones.
Since it is faster to check whether an atom divides an simple element
than to compare two simple elements, we prefer to use~\cite[Algorithm 3]{FG03}
(see Proposition 5.3 there).

\myalg{Alg3}
{Compute $S_\rmright^{[p,q]}(v)$ or $S_\rmleft^{[p,q]}(v)$, respectively.}
{An element $v \in [p,q] \subseteq G^r$, $p,q\in \bbZ^r$.}
{$S_\rmright^{[p,q]}(v)$ (or $S_\rmleft^{[p,q]}(v)$, respectively).}
{\State Let $x_1, \dotsc, x_m$ be the atoms of $G$
\State $R:=\emptyset $
\For{$i=1,\dotsc, m$}
\State Compute $r_{x_i}$ using Algorithm~\ref{rxright} (or~\ref{rxleft}, respectively)
\State $J_i:=\set{ j }{ j\in R \mbox{ and } x_j \preceq r_{x_i} (\mbox{or } r_{x_i} \succeq x_j) }$
\State $K_i:=\set{ j }{ j>i \mbox{ and } x_j \preceq r_{x_i} (\mbox{or } r_{x_i} \succeq x_j) }$
   \If{$J_i=K_i=\emptyset $}
\State $R:=R\cup \{ i\}$
\EndIf
\EndFor
\State \Return{$\set{ r_{x_i} }{ i\in R }$}
}

We conclude by pointing out that, since our invariants are preserved by the
automorphism $\tau$, it is natural to consider them \emph{modulo
$\tau$}, that is, to maintain only one representative (for example, the lexicographically
minimal one) out of each $\tau$-orbit.
For example, the order of $\tau$ for two known Garside structures in braid groups $\BG{N}$
with $N\ge 4$ strands, namely, the Artin--Garside structure with Garside element
$\Delta=\Delta_N$ and the dual or Birman--Ko--Lee structure with $\Delta=\delta_N$, are 2 and $N$,
respectively.

\section{Experimental results in Artin's braid groups}

We have conducted extensive experiments checking the cardinalities of the 
finite sets that can be used for solving the SCP\@.
The experiments are on Artin's braid groups $\BG{N}$, with their
two known Garside structures (Artin and BKL). For a tuple $a\in \BG{N}^r$,
the set $a^{\BG{N}}\cap [\inf a,\infty]$ is the one proposed by 
Lee and Lee~\cite{LL02}. 
The set $a^{\BG{N}}\cap [\inf a,\sup a]$ is its natural subset
introduced here. Both of these sets are not invariants of the conjugacy class.
The sets $\LSS(a)$ and $\LSSS(a)$ are the invariants introduced here,
namely the lexicographic summit set and the lexicographic super summit set.

We did not notice substantial differences between the cardinalities of the two variations
of $\LSSS$ introduced here. Thus, we used in the experiments the second variation,
so that $\LSSS(a)$ is always a subset of $\LSS(a)$. This allows the use of
a smaller number of experiments, while avoiding problems arising
from the large variance.

To give the two sets that are not invariants a fair chance, we considered
them for solving the Search SCP: We constructed conjugate $a,c\in \BG{N}^r$
by choosing $b\in \BG{N}^r$ and $x,y\in \BG{N}$, and setting $a=b^x$ and
$c=b^y$. We then computed, instead of
$a^{\BG{N}}\cap [\inf a,\infty]$, the typically smaller set
$a^{\BG{N}}\cap [\inf c,\infty]$, and similarly for the other set.

Random elements of $\BG{N}$ were generated 
as products of random $2N\log N$ generators, each inverted in
probability $1/2$. Such products are, with high probability, fully
supported in the group.

We summarize the results in Table~\ref{Tab1}, which demonstrates
the following typical inequalities:
\[
\card{\LSSS(a)}
<
\card{\LSS(a)}
\ll
\card{a^{\BG{N}}\cap [\inf c,\sup c]}
<
\card{a^{\BG{N}}\cap [\inf c,\infty]}.
\]
The symbol $\ll$ indicates a dramatic improvement when moving to the invariants.
An additional observation is that the BKL presentation provides much smaller
sets, often one-element sets!

\begin{table}[H]
\label{Tab1}
\caption{Cardinalities of sets (modulo $\tau$) associated to the SCP,
for dimension $r=8$.
Each cell lists the minimum, median, and maximum cardinality encountered,
as well as the percentage of failures,
out of 100 experiments. $\infty$ means $>100{,}000$.}
\begin{tabular}{|r||r|r|r|r|r|r|r|r|}
\cline{2-9}
\multicolumn{1}{c||}{ \bigstrut } &
\multicolumn{2}{c|}{ $a^{\BG{N}}\cap [\inf c,\infty]$ } &
\multicolumn{2}{c|}{ $a^{\BG{N}}\cap [\inf c,\sup c]$ } &
\multicolumn{2}{c|}{ $\LSS(a)$ } &
\multicolumn{2}{c|}{ $\LSSS(a)$ }
\\
\hline
$N$ & Artin & BKL & Artin & BKL & Artin & BKL & Artin & BKL\\
\hline
\hline
4       & 1     & 1 & 1 & 1 & 1 & 1 & 1 & 1\\
       & 116   & 54 & 37 & 17 & 1 & 1 & 1 & 1\\
        & 80{,}438& 27{,}786 & 14{,}318 & 3{,}441 & 8 & 3 & 5 & 1\\
        & 0\% &  0\% &  0\% &  0\% &  0\% &  0\% &  0\% &  0\% \\
\hline
8 & $\infty$ & $\infty$ & 40{,}630 & 872 & 1 & 1 & 1 & 1\\
  & $\infty$ & $\infty$ & $\infty$ & $\infty$ & 63 & 2 & 5   & 1\\
  & $\infty$ & $\infty$ & $\infty$ & $\infty$ & 3{,}732 & 966 & 160 & 17\\
  & 100\% &\ 100\% & 98\% & 96\%   & 0\% & 0\% & 0\% & 0\% \\
\hline
   & $\infty$ & $\infty$ & $\infty$  & $\infty$  & 69{,}534 & 2 & 68 & 1\\
16 & $\infty$   & $\infty$ & $\infty$ & $\infty$ & $\infty$ & 740 & 76{,}509 & 6\\
   & $\infty$ &  $\infty$& $\infty$  & $\infty$  & $\infty$ & $\infty$  & $\infty$ & 28{,}025\\
   & 100\% & 100\% & 100\% & 100\% & 99\% & 10\% & 76\% & 0\%\\
\hline
\end{tabular}
\end{table}

We have tested, for the BKL presentation, the effect of increasing the dimension.
Table 2 summarizes the results. We observe that the cardinality of the invariant $\LSSS$
tends to $1$ with the increase of the dimension, and suggests that when $r=O(N)$
and the elements of the $r$-tuple are ``generic'' and independent, the invariant
tends to have cardinality $1$.

\begin{table}[H]
\caption{The effect of increasing the dimension $r$ on the cardinality of the Birman--Ko--Lee Lexicographic SSS invariant (modulo $\tau$), for braid index $N=32$.
Each cell lists the minimum, median, and maximum cardinality encountered, out of 100 experiments. $\infty$ means $>100{,}000$.}

\begin{tabular}{|l||*{5}{r|}}
\hline
$r$ & 4 & 8 & 16 & 32 & 64\\
\hline
\hline
Minimum & $\infty$         & 720             & 3              & 1 & 1 \\
Median  & $\infty$         & $\infty$        & 95             & 2 & 1 \\
Maximum & $\infty$ & $\infty$ & $\infty$ & 75 & 4 \\
Failures & 100\% & 75\% & 3\% & 0\% & 0\%\\
\hline
\end{tabular}
\end{table}

\section{Comparison with Garside families}
\label{sec:garsidefamilies}

Some of the conditions in
the definition of Garside monoids, like being Noetherian, are often
not needed to establish results about them~\cite{DDM13, DDGM}.
In the monograph on foundations of Garside theory ~\cite{DDGM}, the monoid $M$ is replaced by a left-cancellative 
category $\mathcal{C}$, and the set of simple elements $S$ is generalized to a ``Garside family'', i.e., a subfamily $\mathcal{S}$ of $\mathcal{C}$ such that
every element in $\mathcal{C}$ admits a normal decomposition with respect to $\mathcal{S}$. A \emph{Garside category} is a category $\mathcal{C}$ that admits such a subfamily $\mathcal{S}$.
This generalizes the classic definition in several directions~\cite[Chapter~1]{DDGM}:
\begin{itemize}
	\item[-] The family $\mathcal{S}$ may be infinite (e.g., for the braid group $B_{\infty}$),
	\item[-] the monoid (i.e., the category) may be not Noetherian (e.g., the Klein bottle monoid),
	\item[-] it allows for invertible elements (e.g., wreathed free abelian group $\mathbb{Z}^n \wr S_n$), and
	\item[-] multiplication may not be defined everywhere (e.g., for Ribbon categories).
\end{itemize}
There are several reasons for our working with the classic definition of Garside groups.
Our main focus is the braid group, and for some group theorists and cryptologists, 
the classic definition is more familiar.
We developed formulas for special operations using the
left and right lcm, and established the simultaneous cyclic sliding theorem.
This theorem requires Noetherianity, since it explicitly involves the norm of the Garside element.
To work with Garside categories, we would have to consider cancellative Noetherian Garside categories where every two elements have a left and right lcm, 
and the Garside family is bounded~\cite{DDGM}.
Finally, our algorithms for minimal simple elements require that the number of atoms is finite; we cannot allow infinitely many atoms in $\mathcal{S}$.

The   problem  of  conjugacy  in  the  context  of  Garside  categories  is
treated in chapter VIII of the cited book~\cite{DDGM}. 
Remark~1.16 and Exercise~95 there concern simultaneous conjugacy. 
The solution of the mentioned exercise appears in a work of Digne and Michel~\cite[Proposition~6.2]{DM06}
This establishes the first solution to the simultaneous conjugacy problems for Garside categories, and thus for Garside monoids and groups. 
Restricted to Artin's braid groups, this solution
seems to be considerably less efficient than our solution, since it does not use
minimal intervals or anything similar.

In the context of Garside categories one can also define cycling, decycling and cyclic sliding operations~\cite[Chapter VIII, Definitions~2.3, 2.8, and~2.29]{DDGM}.
These definitions are made in the object category of the conjugacy category. Thus,
a Garside family in the simultaneous conjugacy category does not help defining these operations.
To this end, 
our notions of cycling, decycling and cyclic sliding with respect to a target interval
should be generalized to Garside categories.

\section{Open problems and further work}

In Section~\ref{sec:red} we reduced several problems to the conjunction of
the Search SCP and the computation of the centralizer of a set.
At present, there are no efficient algorithms for the computation of the
centralizers of sets with more than one element in the braid groups.
The computation of the centralizer of an element in
braid groups involves methods used to solve the Conjugacy Problem in these groups~\cite{FG03b}.

\bprb
Does the computation of the centralizer of a set in a group reduce to the
Search SCP?
\eprb

The invariants introduced in the present paper
depend on the order of entries in the $r$-tuple.

\bprb
Is there an invariant, computable in comparable time,
that does not depend on the order of the entries?
\eprb

Our invariants may be huge. In the one-dimensional case, there are
the  much better (essentially, equivalent) invariants of Ultra Summit Sets
and Sliding Circuits. 

\bprb\label{prb:USS}
Is there a generalization of Ultra Summit Sets or Sliding Circuits
to the high-dimensional setting?
\eprb

It is tempting to define, for example, decycling as left conjugation by the left lcm of all tails (or final factors) of all
components.
However, for large simultaneity $r$, the left lcm of all tails is, generically, $\Delta$.
Thus, we conjugate (from the left) by the left lcm of all tails of elements where the supremum can be decreased (i.e., where the element is not in its summit set).
Thus, our decycling is only defined for elements that are outside the (simultaneous) super
summit sets, and they do not define simultaneous USS or sliding circuits in any direct manner.

It is also natural to consider potential applications of this work to cryptanalysis.
As we can see in Section~\ref{sec:red}, the reductions to the Search SCP provide
highly biased instances. The dependency among the entries renders the invariants
too large to be of any direct use. To this end, the invariants must be combined
with heuristic shortcuts, like ones used earlier~\cite{HofSte03}.
The \emph{Search} SCP has the following heuristic speedup: We compute $\LSSS(a)$
and $\LSSS(c)$ in parallel, until we find an element in the intersection.
Heuristically, this has the potential to 
reduce the running time from $n:=\card{\LSSS(a)}$ to
about $\sqrt{n}$. However, in our experiments we did not observe the expected
speedup. An investigation of this phenomenon may help addressing Problem~\ref{prb:USS}.

\ed

\end{document}